\documentclass[a4paper,10pt]{article}

\usepackage{amsmath,enumerate,amsfonts,amssymb,amsthm}

\theoremstyle{plain}
\newtheorem*{Theo}{\bf THEOREM}

\newtheorem{Lem}{\bf LEMMA}

\theoremstyle{definition}
\newtheorem*{Def}{\bf DEFINITION}

\theoremstyle{remark}
\newtheorem*{Rem}{\bf Remark}

\DeclareMathOperator{\re}{Re} \DeclareMathOperator{\im}{Im}
\DeclareMathOperator{\Ker}{Ker}

\begin{document}

\title{\bfseries Potentials for Hyper-K\"ahler Metrics\\
with Torsion} 
\author{Bertrand Banos \and Andrew Swann}
\date{}
\maketitle

\footnotetext[1]{Research supported by the European Human Potential
Programme, HPRN-CT-2000-00101.}  \footnotetext[2]{Banos: Department of
Mathematics and Computer Science, University of Southern Denmark, Campusvej
55, 5230 Odense M, Denmark; email \texttt{banos@imada.sdu.dk}}
\footnotetext[3]{Swann: Department of Mathematics and Computer Science,
University of Southern Denmark, Campusvej 55, 5230 Odense M, Denmark; email
\texttt{swann@imada.sdu.dk}}

\begin{abstract}
  We prove that locally any hyper-K\"ahler metric with torsion admits an
  HKT potential.
\end{abstract}

\section*{Introduction}

A hypercomplex manifold is a manifold endowed with three (integrable)
complex structures $I$, $J$ and $K$ satisfying the quaternion identities
$IJ=-JI=K$.  A metric $g$ compatible with these three complex structures is
said to be hyper-K\"ahler with torsion if the three corresponding K\"ahler
forms satisfy the identities
\begin{equation}{\label{HKT}}
  IdF_I=JdF_J=KdF_K.
\end{equation}
This is equivalent to saying that there exists a connection $\nabla$
preserving the metric ($\nabla g=0$) and the complex structures ($\nabla
I=\nabla J=\nabla K=0$) and whose torsion tensor
\begin{equation*}
  c(X,Y,Z)=g(X,\nabla_YZ-\nabla_ZY-[Y,Z])
\end{equation*}
is totally skew. This connection is necessarily unique and its torsion
tensor is exactly the $3$-form defined by \eqref{HKT}.

The terminology hyper-K\"ahler with torsion is quite misleading since the
underlying metric is in general not K\"ahler at all. We will prefer then
the terminology HKT.

HKT metrics were introduced by Howe and Papadopoulos. They explain in
\cite{HP} how HKT geometry, and other geometries with torsion, arise as the
target spaces of some two-dimensional sigma models in string theory.

Grantcharov and Poon give in \cite{GP} the corresponding mathematical
background. They define in particular the concept of \emph{HKT potential}
which is a natural generalization of the concept of K\"ahler or
hyper-K\"ahler potential. Unlike the K\"ahler case, hyper-K\"ahler metrics
do not admit in general a hyper-K\"ahler potential, even locally (see
\cite{Sw}) but they always admit, locally, an HKT potential. It is actually
easy to check that on a hyper-K\"ahler manifold, any K\"ahler potential for
one of the three complex structures is an HKT potential. We think this
simple remark is sufficient to justify and motivate the question: \emph{do
HKT potentials always exist locally?} We know from Michelson and
Strominger that the answer is yes in the flat case, that is, any HKT metric
on $\mathbb{H}^n$ (with the standard complex structures) admits locally an
HKT potential (see \cite{MS}). It is proved in \cite{PSw} that an HKT
manifold with a special homothety also admits an HKT potential. We prove
here the general result.

Our strategy is based upon the following observation: any compatible metric
on a quaternionic curve (that is, a $4$-dimensional hypercomplex manifold)
is HKT but not necessarily hyper-K\"ahler.  This indicates that HKT
geometry is a better quaternionic generalization of K\"ahler geometry, the
torsion being a direct consequence of the non-commutativity of the
quaternion division ring.  We actually remark that an HKT structure is
essentially a \emph{closed $(1,1)$-form in the sense of Salamon}, that is,
a $2$-form compatible with the three complex structures and closed with
respect to the differential $D$ introduced by Salamon in \cite{Sa}.  This
remark, which already appears in \cite{V} but in other spirit and other
formalism, combined with the properties of the operator $D$ and the twistor
space described by Salamon in \cite{Sa} and \cite{MCS} give us then
directly the wished result.  Indeed the Salamon differential operator used
here provides further analogies with complex geometry and we give
hypercomplex analogues of the local and global \( \partial\bar\partial
\)-lemmas.

\section{HKT metrics and Salamon $(1,1)$-forms}

An almost hypercomplex manifold is a (smooth) manifold endowed with $3$
almost complex structures $I$, $J$ and $K$ satisfying the quaternion
identities
\begin{equation*}
  IJ=-JI=K.
\end{equation*}
Note that on a almost hypercomplex manifold there is actually a $2$-sphere
worth of almost complex structures:
\begin{equation*}
  S^2=\{aI+bJ+cK: a^2+b^2+c^2=1\}.
\end{equation*}
The integrability of $I$, $J$ and $K$ is equivalent to the existence of a
(unique) torsion-free connection $\nabla^{\text{Ob}}$ preserving the
quaternion action, the so-called Obata connection.

\subsection{Exterior forms}

Let $(M,I,J,K)$ be an almost hypercomplex manifold and let $\Lambda^k$ be
the bundle of $k$-forms on $M$. We denote by $\Lambda^{p,q}_\mathcal{I}$
the subbundle of forms of type $(p,q)$ with respect to the almost complex
structure $\mathcal{I}\in S^2$.

Studying the action of $GL(n,\mathbb{H})$ on $\Lambda^k(T^*)$, Salamon
introduces in \cite{Sa} the subbundle
\begin{equation*}
  A^k=\underset{\mathcal{I}\in
  S^2}{\sum}\big(\Lambda^{k,0}_\mathcal{I}\oplus
  \Lambda^{0,k}_\mathcal{I}\big).
\end{equation*}
This bundle can be understood as the analogue for hypercomplex manifolds of
the bundle $\Lambda^{k,0}\oplus \Lambda^{0,k}$ for complex manifolds (see
\cite{W}).

It will be convenient for us to choose a preferred complex structure, say
$I$. Although this choice is not really natural (all the complex structures
have the same status and should be studied together), we will see that it
is useful for understanding HKT geometry as a quaternionic analogue of
K\"ahler geometry. When $M$ is considered as a complex manifold, \emph{this
is always with respect to the complex structure $I$.} We will write for
example $\Lambda^{p,q}$ for $\Lambda^{p,q}_I$. The Hodge decomposition of
$\Lambda^2$ with respect to $I$ induces the decomposition
\begin{equation*}
  A^2=\Lambda^{2,0}\oplus \Lambda^{0,2}\oplus A^{1,1},
\end{equation*}
with
\begin{equation*}
  A^{1,1}=\big\{\omega\in \Lambda^2:\text{ $I\omega=\omega$ and
  $J\omega=-\omega$}\big\}.
\end{equation*}
For example, if $g$ is a hyperhermitian metric on $M$ then the K\"ahler
form $F_I=g(I\cdot,\cdot)$ is a smooth section of $A^{1,1}$ and conversely
any smooth section $F_I$ of $A^{1,1}$ defines an (possibly indefinite
and/or degenerate) hyperhermitian metric $g=-F_I(I\cdot,\cdot)$. We will
call such a form a $(1,1)$-form in the sense of Salamon.

\subsection{The Salamon differential}

There is an orthogonal projection $\eta\colon\Lambda^k\rightarrow A^k$
whose kernel is the subbundle
\begin{equation*}
  B^k=\bigcap_{\mathcal{I}\in S^2}
  \big(
  \Lambda^{k-1,1}_\mathcal{I} \oplus
  \Lambda^{k-2,2}_\mathcal{I} \oplus \dots \oplus
  \Lambda^{1,k-1}_\mathcal{I}
  \big).
\end{equation*}
Let $\mathcal{A}$ denote the space of smooth sections of $A$. The
Salamon differential
\begin{equation*}
  D\colon\mathcal{A}^k\rightarrow\mathcal{A}^{k+1}
\end{equation*}
is simply the composition of the projection $\eta$ with the de Rham
differential $d$:
\begin{equation*}
  D=\eta\circ d.
\end{equation*}
For example, if $\theta$ is a $1$-form on $M$ then
\begin{equation}{\label{proj}}
  D\theta= (d\theta)^{2,0} + (d\theta)^{0,2} +
  \tfrac12\big((d\theta)^{1,1} - J(d\theta)^{1,1}\big).
\end{equation}

Salamon shows in \cite{Sa} the following:

\begin{Theo}[Salamon]
  An almost hypercomplex structure is integrable if and only if $D^2=0$.
\end{Theo}

This result is completely analogous to the corresponding statement
involving an almost complex structure and the Dolbeault operator
$\overline{\partial}$.

\subsection{The twistor space}

If $(M,I,J,K)$ is a hypercomplex manifold then the manifold $Z=M\times S^2$
admits an integrable complex structure $\mathbb{I}$ defined by
\begin{equation*}
  \mathbb{I}_{(x,\vec{a})}=
  \begin{pmatrix}
    a_1 I_x + a_2 J_x + a_3 K_x & 0 \\
    0 & I_{\vec{a}} \\
  \end{pmatrix}
  ,
\end{equation*}
where $I_{\vec{a}}\colon X\mapsto \vec{a}\times X$ is the usual complex
structure on $T_{\vec{a}}S^2$. The space $Z$ endowed with the complex
structure $\mathbb{I}$ is called the twistor space of the hypercomplex
manifold $M$.

The cohomology of the twistor space can be related to the cohomology of the
Salamon's complex as follows (\cite{MS}):

\begin{Theo}[Mamone Capria, Salamon]
  Let $M^{4n}$ be a hypercomplex manifold with twistor space $Z$.
  Then
  \begin{equation*}
    H^k(Z,\mathcal{O})\cong
    \begin{cases}
      \frac{\Ker(D\colon\mathcal{A}^k\rightarrow
      \mathcal{A}^{k+1})}{\im(D\colon\mathcal{A}^{k-1}\rightarrow
      \mathcal{A}^{k})},& 0\leqslant k\leqslant 2n,\\
      0,& k=2n+1.
    \end{cases}
  \end{equation*}
\end{Theo}

\subsection{HKT metrics}

Let $(M,I,J,K)$ be a hypercomplex manifold and let $g$ be a hyperhermitian
metric on $M$, that is
\begin{equation*}
  g(IX,IY)=g(JX,JY)=g(KX,KY)=g(X,Y),
\end{equation*}
for all tangent vectors $X$ and $Y$. We will denote by $F_\mathcal{I}$ the
K\"ahler form associated with the complex structure $\mathcal{I}$:
\begin{equation*}
  F_\mathcal{I}=g(\mathcal{I}\cdot,\cdot).
\end{equation*}
Note that $g$ can be indefinite in what follows.

In the physics literature, this hyperhermitian metric is said to be HKT if
there exists a hyperhermitian connection whose torsion tensor is totally
antisymmetric (see \cite{HP}). We will rather use the reformulation
introduced by Grantcharov and Poon:

\begin{Def}{\label{Def}}
  The hyperhermitian metric $g$ is \emph{HKT} if
  \begin{equation*}
    IdF_I=JdF_J=KdF_K.
  \end{equation*}
\end{Def}

For example, any compatible metric on a quaternionic curve is HKT:

\begin{Lem}
  Any hyperhermitian metric on a $4$-dimensional hypercomplex manifold is
  HKT.
\end{Lem}

\begin{proof}
  Let $(M,I,J,K)$ be a $4$-dimensional manifold and let $g$ be a
  hyperhermitian metric on it. It is noted in \cite{PS} that $g$ is
  necessarily Einstein-Weyl with respect to the Obata connection
  $\nabla^{\text{Ob}}$.  In particular, there exists a $1$-form $\omega$
  such that
  \begin{equation*}
    \nabla^{\text{Ob}} g =\omega\otimes g.
  \end{equation*}
  (This may be seen directly, by noting that on a four-manifold the
  conformal class of~$g$ is uniquely determined by $I$, $J$ and~$K$.)
  Since $\nabla^{\text{Ob}}$ is torsion-free and compatible with $I$, $J$
  and $K$ this implies that
  \begin{equation*}
    dF_I=\omega\wedge F_I, \quad dF_J=\omega\wedge F_J, \quad
    dF_K=\omega\wedge F_K.
  \end{equation*}
  Define now $\alpha$ by $\omega=\lambda\alpha$ and $\|\alpha\|=1$.  Since
  $M$ is $4$-dimensional we get
  \begin{equation*}
    \left\{
      \begin{aligned}
        F_I&=\alpha\wedge I\alpha + J\alpha\wedge K\alpha,\\
        F_J&=\alpha\wedge J\alpha + K\alpha\wedge I \alpha,\\
        F_K&=\alpha\wedge K\alpha + I\alpha\wedge J\alpha.
      \end{aligned}
    \right.
  \end{equation*}
  And then
  \begin{equation*}
    IdF_I=JdF_J=KdF_K=\lambda I\alpha\wedge J\alpha\wedge K\alpha.
  \end{equation*}
\end{proof}

Many of the explicitly known HKT examples in higher dimensions are
homogeneous and come from the Joyce hypercomplex structures associated to
any compact semi-simple Lie group (see \cite{J} and \cite{GP}). For
example, the Killing-Cartan metric on $SU(3)$ is HKT for the (non-trivial)
invariant hypercomplex structure on $SU(3)$ constructed by Joyce.  It is
worth mentioning that the Lie bracket on $su(3)$ is exactly the torsion of
the HKT structure. In particular, due to the Jacobi identity, the torsion
form is closed: $SU(3)$ is a strong HKT manifold.

\subsection{HKT forms}

When one is more interested in complex and symplectic properties than in
Riemannian ones, one can define a K\"ahler structure as a non-degenerate
closed $(1,1)$-form. It is possible to have a similar approach for HKT
structures.  The following result is due to Verbitsky~\cite{V}, but we
prefer to give a direct proof using the Obata connection.

\begin{Lem}{\label{keylemma}}
  Let $F\in \mathcal{A}^{1,1}$ be a non-degenerate Salamon (1,1)-form on a
  hypercomplex manifold $(M,I,J,K)$ . The (pseudo) metric
  \begin{equation*}
    g=-F(I\cdot,\cdot)
  \end{equation*}
  is HKT if and only if $F$ is $D$-closed:
  \begin{equation*}
    DF=0.
  \end{equation*}
  Such a form is called an HKT form.
\end{Lem}

\begin{proof}
  Suppose that $g$ is HKT. For any complex structure $\mathcal{I}\in S^2$
  the form $dF_\mathcal{I}$ has type $(2,1)+(1,2)$ with respect to the
  complex structure $\mathcal{I}$. But since $IdF_I=JdF_J=KdF_K$ we deduce
  that $dF_I$ has type $(2,1)+(1,2)$ with respect to the three complex
  structures: $dF_I\in \mathcal{B}^3$ that is $DF_I=0$. Since $F_I=F$, we
  obtain the result.

  Suppose now that $DF=0$. This is equivalent to the relation
  \begin{equation*}
    dF(U,V,W)=dF(\mathcal{I}U,\mathcal{I}V,W) +
    dF(\mathcal{I}U,V,\mathcal{I}W) + dF(U,\mathcal{I}V,\mathcal{I}W),
  \end{equation*}
  for all $\mathcal{I}$ in $S^2$. In particular
  \begin{equation}{\label{eq1}}
    dF(IU,IV,IW)=dF(KU,KV,IW)+ dF(KU,IV,KW)+dF(IU,KV,KW).
  \end{equation}
  Since the Obata connection is torsion-free, the following holds:
  \begin{equation}{\label{eq2}}
    dF(X,Y,Z)=\nabla^{\text{Ob}} F(X,Y,Z)+ \nabla^{\text{Ob}} F(Y,Z,X) +
    \nabla^{\text{Ob}} F(Z,X,Y).
  \end{equation}
  Moreover, since $F\in \mathcal{A}^{1,1}$, we have
  \begin{equation}{\label{eq3}}
    \left\{
      \begin{aligned}
        &F(X,Y)=F(IX,IY)=-F(KX,KY),\\
        &F(IX,KY)=F(KX,IY).
      \end{aligned}
    \right.
  \end{equation}
  Using $\eqref{eq2}$ and $\eqref{eq3}$ in $\eqref{eq1}$ we obtain
  \begin{equation*}
    \begin{split}
      &dF(IU, IV,IW)\\
      &=2\nabla^{\text{Ob}} F(KU,KV,IW) + 2\nabla^{\text{Ob}}
      F(KV,KW,IU) + 2\nabla^{\text{Ob}} F(KW,KU,IV)\\
      &\qquad - \nabla^{\text{Ob}} F(IU,V,W) - \nabla^{\text{Ob}} F(IV,W,U)
      - \nabla^{\text{Ob}} F(IW,U,V)
    \end{split}
  \end{equation*}
  and thus
  \begin{equation*}
    \begin{split}
      dF(IU,IV,IW)
      &=\nabla^{\text{Ob}} F(KU,KV,IW)+\nabla^{\text{Ob}} F(KV,KW,IU)\\
      &\qquad+ \nabla^{\text{Ob}} F(KW,KU,IV).\\
    \end{split}
  \end{equation*}
  Define now $G=F(J\cdot,\cdot)$.
  \begin{equation*}
    \begin{split}
      &dG(KU, KV,KW)\\
      &\quad = \nabla^{\text{Ob}} G(KU,KV,KW) + \nabla^{\text{Ob}}
      G(KV,KW,KU) + \nabla^{\text{Ob}} G(KW,KU,KV)\\
      &\quad = \nabla^{\text{Ob}} F(KU,IV,KW) + \nabla^{\text{Ob}}
      F(KV,IW,KU) + \nabla^{\text{Ob}} F(KW,IU,KV)\\
      &\quad = \nabla^{\text{Ob}} F(KU,KV,IW) + \nabla^{\text{Ob}}
      F(KV,KW,IU) + \nabla^{\text{Ob}} F(KW,KU,IV)
    \end{split}
  \end{equation*}
  and therefore $dF(IU,IV,IW)=dG(KU,KV,KW)$. In other words, $IdF_I=KdF_K$
  with $F_I=F$ and $F_K=G$.
\end{proof}

\begin{Rem}
  We know from Fino and Grantcharov (\cite{FG}) that there exists some
  hypercomplex manifolds which do not admit an HKT metric. Lemma
  \ref{keylemma} seems to indicate that the question of existence of an HKT
  metric on a given hypercomplex manifold is highly non-trivial.
\end{Rem}

\section{HKT potentials}

Let $(M,I,J,K)$ be a hypercomplex manifold. Following $\cite{GP}$ we define
the action of $\mathcal{I}\in S^2$ on $k$-forms by
\begin{equation*}
  \mathcal{I}\omega(X_1,\dots,X_k) =
  (-1)^k\omega(\mathcal{I}X_1,\dots,\mathcal{I}X_k)
\end{equation*}
and the differential $d_k$ is
\begin{equation*}
  d_\mathcal{I}\omega=(-1)^k\mathcal{I}d\mathcal{I}\omega.
\end{equation*}
Note that $d$, $d_I$, $d_J$ and $d_K$ all anti-commute.

Recall that a hyperhermitian metric $g$ on $M$ is said to be hyper-K\"ahler
if it is K\"ahler for each complex structure. A possibly locally defined
function $\mu$ is a hyper-K\"ahler potential for this metric $g$ if it is a
K\"ahler potential for each complex structure, that is,
\begin{equation*}
  F_I=dd_I\mu,\quad F_J=dd_J\mu,\quad F_K=dd_K\mu.
\end{equation*}
It is proved in \cite{Sw} that such a potential does not exist in general
but it is straightforward to check that if $\nu$ is a K\"ahler potential
for the complex structure $I$ then
\begin{equation*}
  F_I=dd_I\nu,\quad F_J=\tfrac12(dd_J+d_Kd_I)\nu,\quad F_K=\tfrac12(dd_K+
  d_Id_J)\nu.
\end{equation*}
We say then that any hyper-K\"ahler metric admits an HKT potential:

\begin{Def}[Grantcharov, Poon] A possibly locally defined function $\mu$ is
  an \emph{HKT potential} for an HKT metric $g$ if
  \begin{equation*}
    F_I=\tfrac12(dd_I+d_Jd_K)\mu,\quad
    F_J=\tfrac12(dd_J+d_Kd_I)\mu,\quad
    F_K=\tfrac12(dd_K+d_Id_J)\mu.
\end{equation*}
\end{Def}

\begin{Rem}
  Note that on an HKT manifold the following identities are actually
  equivalent:
  \begin{enumerate}
  \item $F_I=\frac12(dd_I+d_Jd_K)\mu$,
  \item $F_J=\frac12(dd_J+d_Kd_I)\mu$,
  \item $F_K=\frac12(dd_K+d_Id_J)\mu$,
  \item $g=\frac14(1+I+J+K)(\nabla^{\text{Ob}})^2\mu$.
  \end{enumerate}
\end{Rem}

\subsection{The four-dimensional case}

In this dimension one can check directly that HKT metrics always admit an
HKT potential:

\begin{Lem}
  Let $g$ be an HKT metric on a $4$-dimensional hypercomplex manifold and
  let $\omega$ be the $1$-form defined by the Obata connection via
  $\nabla^{\textup{Ob}}g=\omega\otimes g$ .

  A function $\mu$ is an HKT potential for $g$ if and only if it is
  solution of the elliptic equation
  \begin{equation*}
    \Delta\mu-\omega^{\sharp}(\mu)+4=0,
  \end{equation*}
  where \( \Delta \) is the Laplacian of the Riemannian metric~$g$.
\end{Lem}

Local existence of HKT potentials now follows from the general theory for
the Laplace operator, see for example \cite{GT}.

\begin{proof}
  Let $\nabla^{\text{LC}}$ be the Levi-Civita connection and define
  $a=\nabla^{\text{Ob}}-\nabla^{\text{LC}}$. Since
  $\nabla^{\text{Ob}}g=\omega\otimes g$ and $\nabla^{\text{LC}}g=0$ we get
  \begin{equation*}
    \omega(U)g(V,W)=-g(a_UV,W)-g(a_UW,V),
  \end{equation*}
  for all vector fields $U$, $V$ and $W$. Moreover, $a_UV=a_VU$ holds for
  all $U$ and $V$ since $\nabla^{\text{Ob}}$ and $\nabla^{\text{LC}}$ are
  torsion free. We now obtain
  \begin{equation*}
    g(a_UV,W)=\tfrac12\big( \omega(W)g(U,V) - \omega(U)g(V,W) -
    \omega(V)g(U,W) \big),
  \end{equation*}
  for all $U$, $V$ and $W$. In particular, if $X$ is a (local) unit vector
  field, then
  \begin{equation*}
    g(a_XX,Y)=\tfrac12\omega(Z)-\omega(X)g(X,Z)
  \end{equation*}
  and
  \begin{equation*}
    g(a_XX+a_{IX}IX+ a_{JX}JX+a_{KX}{KX},Z)=\omega(Z)
  \end{equation*}
  for all $Z$.

  The metric $g$ is the unique hyperhermitian metric satisfying $g(X,X)=1$.
  Therefore $\mu$ is an HKT potential if and only if
  \begin{equation*}
    \tfrac14(1+I+J+K)(\nabla^{\text{Ob}})^2\mu(X,X)=1,
  \end{equation*}
  that is,
  \begin{equation*}
    \text{Trace}(\nabla^{\text{Ob}}d\mu)=4.
  \end{equation*}
  Note that the Laplacian $\Delta\mu$ is by definition
  $-\text{Trace}(\nabla^{\text{LC}}d\mu)$. Thus $\mu$ is a HKT potential
  for $g$ if and only if
  \begin{equation*}
    -\Delta \mu + d\mu(a_XX+a_{IX}IX+a_{JX}JX+a_{KX}{KX})=4.
  \end{equation*}
\end{proof}

\subsection{The local $DD_I$-lemma}

The easiest way to show that a K\"ahler metric admits a local K\"ahler
potential is to apply the local $dd_I$-lemma to the closed (and therefore
locally exact) K\"ahler form. This is exactly the same for HKT potentials
if one now uses the Salamon differential:

\begin{Lem}
  A HKT metric locally admits a potential if and only if the corresponding
  HKT form is locally $D$-exact.
\end{Lem}

\begin{proof}
  Suppose that $F=\frac12(dd_I+d_Jd_K)\mu$. Then
  $F=\frac12(d\theta-Jd\theta)$ with $\theta=Id\mu$. Note that $d\theta$ is
  a $(1,1)$-form (for $I$) since $d\theta=dd_I\mu$.  Therefore, according
  to \eqref{proj} $F=D\theta$.

  Conversely, suppose that $F=D\theta$ for some $1$-form $\theta$.  Since
  $F$ is a $(1,1)$-form for $I$, we obtain from \eqref{proj}
  \begin{equation*}
    \left\{
      \begin{aligned}
        &d\theta\in \Lambda^{1,1}, \\
        &F=\tfrac12(d\theta-Jd\theta).
      \end{aligned}
    \right.
  \end{equation*}
  Since $I$ is an integrable complex structure, the local $dd_I$-lemma
  holds: locally there exists $\mu$ such that $d\theta=dd_I\mu$. We get
  then
  \begin{equation*}
    F=\tfrac12(dd_I-Jdd_I)\mu=\tfrac12(dd_I+d_Jd_K)\mu.
  \end{equation*}
\end{proof}

\begin{Theo}
  Any HKT metric admits locally an HKT potential.
\end{Theo}

\begin{proof}
  Let $g$ be an HKT metric on a hypercomplex manifold $(M,I,J,K)$ and let
  $F=F_I$ be the corresponding HKT form.  This form is $D$-closed and
  according to the theorem of Mamone Capria and Salamon it implies that it
  is locally $D$-exact. The idea of the proof is the following:

  Let $Z=M\times S^2$ be the twistor space of $M$ and $p\colon Z\rightarrow
  M$ the natural projection. Define the $(0,2)$-form $G$ on $Z$ by
  \begin{equation*}
    G_{(x,\vec{a})}=\big(F_x\big)^{0,2}_{\vec{a}}.
  \end{equation*}
  The form $G$ in the point $(x,\vec{a})$ is the $(0,2)$-part of the form
  $F$ in the point $x$ with respect to the complex structure
  $\vec{a}=a_1I_x+a_2J_x+a_3K_x$. Grantcharov and Poon have proved in
  \cite{GP} that $g$ is HKT if and only if the form
  $\big(F\big)^{0,2}_{\vec{a}}$ is a $\overline{\partial}_{\vec{a}}$-closed
  form on $M$. Moreover $\big(F\big)^{0,2}_{\vec{a}}$ is holomorphic in
  $\vec{a}$. This implies that $G$ is a $\overline{\partial}$-closed
  $(0,2)$-form on $Z$. Now a 1-pseudo-convexity argument says that one can
  always choose a neighbourhood $U$ of a point $x\in M$ such that
  \begin{equation*}
    H^{0,2}_{\overline{\partial}}(p^{-1}(U))=H^2(p^{-1}(U),\mathcal{O})=0.
  \end{equation*}
  It implies that it exists a $1$-form $\phi$ on $p^{-1}U$ such that
  $G=\overline{\partial}\phi$. Moreover one can choose this form without
  part on $S^2$: for any point $x\in U$, $\phi_x$ is a holomorphic section
  of the bundle over $S^2$ with fibre
  \begin{equation*}
    \mathcal{B}_{\vec{a}} = \big\{ \omega+i\vec{a}\omega :
    \omega\in T_x^*M \big\}.
  \end{equation*}
  Using now the compactness of $S^2$ we deduce that, in any point $\vec{a}$
  of $S^2$, $\phi=\theta+i\vec{a}\theta$ with $\theta$ a $1$-form on $U$.
  We get then $\re(G_{\vec{a}})=\frac12(d\theta-\vec{a}\theta)$ for any
  $\vec{a}\in S^2$. Taking $\vec{a}=I$ and $\vec{a}=J$ we obtain
  \begin{equation*}
    \left\{
      \begin{aligned}
        &d\theta-Id\theta=0,\\
        &F=\tfrac12(d\theta-Jd\theta),\\
      \end{aligned}
    \right.
  \end{equation*}
  that is, $F=D\theta$ on the neighbourhood $U$ of a fixed point $x$.
\end{proof}

\begin{Rem}
  This actually shows that the local $DD_I$-lemma holds on hypercomplex
  manifolds with $D_I=(-1)^k IDI$.
\end{Rem}

\subsection{The global $DD_I$-lemma}

Let $(M,g,I,J,K)$ be an HKT manifold with HKT form $F$. As $F$ is
$D$-closed it defines a Salamon cohomology class $[F]\in H^2_D(M)$ which we
can call the \emph{HKT class}. Assume that $F'$ is another HKT form in the
same HKT class, that is $F-F'=D\theta$.  Since $F-F'$ is a $(1,1)$-form
(for $I$), we get
\begin{equation*}
  F-F'=\tfrac12(d\theta-Jd\theta),
\end{equation*}
with $d\theta\in \Lambda^{1,1}$. Therefore, \emph{if the global
$dd_I$-lemma holds on $M$} then there exists a global function $\phi$ on
$M$ such that $F'=F+DD_I\phi$.

Note that if $DD_I\phi=0$ then $\phi$ is harmonic with respect the complex
Laplacian $\Delta^c$ defined by
\begin{equation*}
  \overline{\partial}^\star\overline{\partial} f =\Delta^c f = g(dd_If,F_I).
\end{equation*}
Indeed if $DD_If=0$ then $dd_If= Jdd_If$ and then
\begin{equation*}
  \Delta^c f = g(dd_If,F_I) = g(Jdd_I,F_I) = g(dd_I,JF_I) =
  -g(dd_I\phi,F_I) = -\Delta^c f.
\end{equation*}

Finally we get the following:

\begin{Theo}
  Let $(M,I,J,K)$ be a compact hypercomplex manifold on which the global
  $dd_I$-lemma holds and let $g$ and $g'$ be two HKT metrics with same HKT
  class $[F]=[F']$. Then there exists a smooth real function $\phi$ on $M$
  such that $F'=F+DD_I\phi$.  This function is unique up to a constant.
  \qed
\end{Theo}


\begin{thebibliography}{99}

\bibitem[FG]{FG} A. Fino, G. Grantcharov: \emph{On some properties of the
  manifolds with skew-symmetric torsion and holonomy SU(n) and Sp(n)},
  preprint math.DG/0302358, in corso di stampa su Advances in Mathematics

\bibitem[GP]{GP} G. Grantcharov, Y.S. Poon: \emph{Geometry of
  hyper-K\"ahler connections with torsion}, Comm. Math. Phys. 213 (2000),
  No.~1, p. 19--37

\bibitem[GT]{GT} D. Gilbarg, N.S. Trudinger: \emph{Elliptic partial
  differential equations of second order}, Grundlehren der Mathematischen
  Wissenschaften vol.~224, Second Edition, Springer-Verlag, Berlin, 1983
  
\bibitem[HP]{HP} P.S. Howe, G. Papadopoulos: \emph{Twistor spaces for
  hyper-K\"ahler manifolds with torsion}, Phys. Lett., B 379 (1996),
  No.~1-4, p. 80--86

\bibitem[J]{J} D. Joyce: \emph{Compact hypercomplex and quaternionic
  manifolds}, J. Differential Geom. 35 (1992), p.743--761

\bibitem[MCS]{MCS} M. Mamone Capria, S.M. Salamon: \emph{Yang-Mills on
quaternionic spaces}, Nonlinearity 1 (1988), p. 517--530

\bibitem[MS]{MS} J. Michelson, A. Strominger: \emph{The geometry of (super)
  conformal quantum mechanics}, Commun. Mat. Phys. 213 (2000), p. 1--17

\bibitem[PS]{PS} H. Pedersen, A. Swann: \emph{Riemannian submersions,
  four-manifolds and Einstein-Weyl geometry}, Proc.  London Math. Soc. 66
  (1993), p. 381--399

\bibitem[PSw]{PSw} Y.S. Poon, A. Swann: \emph{Potential functions of HKT
  spaces}, Classical and Quantum Gravity 18 (2001), p.  4711--4714

\bibitem [Sa]{Sa} S.M. Salamon: \emph{Differential geometry
of quaternionic manifolds}, Ann. Sc. Ec. Norm. Sup $4^e$ s\'erie,
19 (1986), p.31--55

\bibitem [Sw]{Sw} A. Swann: \emph{HyperK\"ahler and quaternionic K\"ahler
  geometry}, Math. Ann. 289 (1991), p.  421--450

\bibitem [V]{V} M. Verbitsky: \emph{HyperK\"ahler manifolds with torsion,
  Supersymmetry and Hodge theory}, Asian J. Math. vol. 6 (2002), p.
  679--712

\bibitem[W]{W} D. Widdows: \emph{A Dolbeault-type double complex on
  quaternionic manifolds}, Asian J. Math., vol. 6, No.~2 (2002), p.
  253--276
\end{thebibliography}
\end{document}